\documentclass[12pt]{article}

\usepackage{amsfonts, amsmath, amssymb}
\usepackage{graphicx} 
\usepackage{float} 
\usepackage{xcolor,pict2e}
\usepackage{enumitem}
\usepackage[cbgreek]{textgreek}
\usepackage{hyperref}

\global\arraycolsep=2pt
\allowdisplaybreaks

\newcommand\he{{\hat{e}}}

\newcommand\lalpha{{\underline{\alpha}}}
\newcommand\lbeta{{\underline{\beta}}}

\newcommand\lsigma{{\underline{\sigma}}}
\newcommand\lmu{{\underline{\mu}}}
\newcommand\lnu{{\underline{\nu}}}
\newcommand{\lx}[1]{{\underline{#1}}}

\newcommand\tI{{\mathtt{I}}}
\newcommand\tJ{{\mathtt{J}}}
\newcommand\tK{{\mathtt{K}}}
\newcommand\tL{{\mathtt{L}}}
\newcommand\talpha{{\textrm{\textalpha}}}
\newcommand\tbeta{{\textrm{\textbeta}}}
\newcommand\tgamma{{\textrm{\textgamma}}}
\newcommand\tsigma{{\textrm{\textsigma}}}
\newcommand\tlambda{{\textrm{\textlambda}}}
\newcommand\tdelta{{\textrm{\textdelta}}}
\newcommand\tmu{{\textrm{\textmu}}}
\newcommand\tnu{{\textrm{\textnu}}}
\newcommand\ti{{\mathtt{i}}}
\newcommand\tj{{\mathtt{j}}}
\newcommand\tk{{\mathtt{k}}}
\newcommand\tl{{\mathtt{l}}}
\newcommand\tm{{\mathtt{m}}}
\newcommand{\tx}[1]{{\mathtt{#1}}}

\newcommand\bI{{\boldsymbol{I}}}
\newcommand\bJ{{\boldsymbol{J}}}
\newcommand\bK{{\boldsymbol{K}}}
\newcommand\balpha{{\boldsymbol{\alpha}}}
\newcommand\bbeta{{\boldsymbol{\beta}}}

\newcommand\bsigma{{\boldsymbol{\sigma}}}
\newcommand\bmu{{\boldsymbol{\mu}}}
\newcommand\bnu{{\boldsymbol{\nu}}}
\newcommand\bi{{\boldsymbol{i}}}
\newcommand\bj{{\boldsymbol{j}}}
\newcommand\bk{{\boldsymbol{k}}}
\newcommand{\bx}[1]{{\boldsymbol{#1}}}

\newcommand\be{{\mathbf{e}}}

\begin{document}

\title{\fontsize{18}{24} \bf Moving frame and spin field representations of submanifolds in flat space}
\author{\normalsize Shou-Jyun Zou\thanks{sgzou2000@gmail.com}}
\date{\today}

\maketitle

\begin{abstract}
We introduce a spin field approach, that is compatible with the Cartan moving frame method, to describe the submanifold in a flat space.
In fact, we consider a kind of spin field $\psi$, that satisfies a Killing spin field equation (analogous to a Killing spinor equation) written in terms of the Clifford algebra, and we use the spin field to locally rotate the orthonormal basis $\{\hat{e}_\mathtt{I}\}$.
Then, the deformed orthonormal frame $\{\tilde{\psi}\hat{e}_\mathtt{I}\psi\}$ can be seen as the moving frame of a submanifold.
We find some solutions to the Killing spin field equation and demonstrate an explicit example.
Using the product of the spin fields, one can easily generate a new immersion submanifold, and this technique should be useful for studies in geometry and physics.
Through the spin field, we find a linear relation between the connection and the extrinsic curvature of the submanifold.
We propose a conjecture that any solution of the Killing spin field equation can be locally written as the product of the solutions we find.
\end{abstract}

\tableofcontents

\section{Introduction}

The isometric embedding of a (pseudo-)Riemannian manifold into an ambient flat space is a long-standing problem.
There are several famous results, such as the Janet-Cartan theorem \cite{Janet:1926,Cartan:1927} for locally isometric embedding (immersion) problem, and the Nash theorem \cite{Nash:1956} for gobally isometric  embedding problem.
These results show that we can describe any Riemannian manifold in the higher enough dimensional flat space.
This provides some advantages and more intuitive ways of solving the mathematical and physical problems.

In this article, we use the moving frame \cite{Cartan:2001} and the spin field\footnotemark[2] approaches to describe the 4-dimensional submanifold in a 10-dimensional flat space, due to the Janet-Cartan theorem, that can locally describe all 4-dimensional manifolds.
We introduce a kind of spin field $\psi$, that satisfies a Killing spin field equation (analogous to a Killing spinor equation) written in terms of the Clifford algebra, and we use the spin field to locally rotate the orthonormal basis $\{\he_\tI\}$.
Then the deformed orthonormal frame $\{\tilde{\psi}\he_\tI\psi\}$ can be seen as the moving frame of a submanifold in a flat space.
\footnotetext[2]{
    Our spin field approach is inspired by the spinor representations of submanifolds \cite{Vaz:2019umf,Bayard:2017}.
    Usually, the definition of a spinor is the minimal left ideal of the spin group.
    Our approach only requests the field to have the properties of the spin group, so we call it the spin field.
}

We find some solutions to the Killing spin field equation and demonstrate an explicit example.
Using the product of the spin fields, we can easily generate a new immersion submanifold, and this technique should be useful for studies in geometry and physics \cite{Willison:2013ova}.
Through the spin field, we find a linear relation between the connection and the extrinsic curvature of the submanifold.
A conjecture is proposed that any solution of the Killing spin field equation can be locally written as the product of the solutions we find.

\section{Isometric  immersion}

To provide a clear description, here we focus on a pseudo-Riemannian manifold $\mathcal{M}^4$ with signature $(-,+,+,+)$ isometric  immersion in a 10-dimensional Minkowski space $\mathbb{R}^{1,9}$, and the method can be easily generalized to any dimensional case.
The isometric immersion map is denoted by $\bx{q}: \mathcal{M}^4\to \mathbb{R}^{1,9}$.

In this paper, we will use 3 kinds of coordinate systems to describe the geometry $\mathcal{M}^4$ in $\mathbb{R}^{1,9}$, and we use different fonts to display their indices: 
\begin{description}
    \item [The underline font] ($x^{\lmu}$) for the local coordinate of the immersion manifold : \\
    The underline font ``$\lmu$'' denotes the index of the local coordinate of the immersion submanifold,  e.g. $x^{\lmu}$, $x^{\underline{0}}$. 
    The range of the index $\lmu$ is $\{0,1,2,3\}$.
    \item [The normal text font] ($\he_{\tI}$) for the orthonormal fixed frame : \\
    The normal text font ``$\tI$'' denotes the index of the orthonormal fixed frame, e.g. $\he_{\tI}$, $\he_{\mathtt{1}}$.
    The range of the index $\tI$ is $\{0,1,\dots,9\}$.
    \item [The bold font] ($\be_{\bI}$) for the orthonormal moving frame : \\
    The bold font ``$\bI$'' denotes the index of the orthonormal moving frame, e.g. $\be_{\bI}$, $\be_{\boldsymbol{1}}$.
    The range of the index $\bI$ is $\{0,1,\dots,9\}$.
\end{description}
The following schematic diagram shows a submanifold with the local coordinates, the orthonormal fixed frame and the orthonormal moving frame:
\begin{figure}[H]
    \centering
    \includegraphics[width=0.6\textwidth]{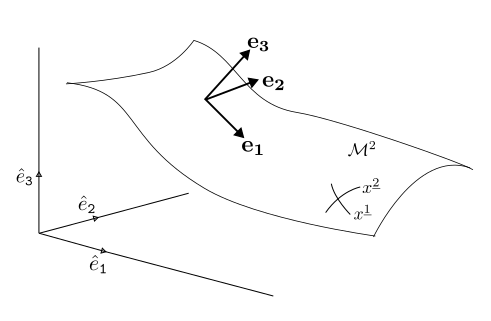}
    \caption{The example for a submanifold with different coordinate systems.}
\end{figure}

In this article, we assume that the first four basis vectors $\be_\bmu$ ($\bmu\in\{0,1,2,3\}$) of the moving frame are tangent to the submanifold, and the remaining basis vectors $\be_\bi$ ($\bi\in\{4,5,\dots,9\}$) are normal to the submanifold.
To provide a clearer description of the space related to submanifolds, we use three types of characters to represent and specify the indices belonging to different ranges:
\begin{description}
    \item [Uppercase Latin characters] $I\in\{ 0,1,\dots,9 \}$ : \\
    All indices of the basis of the ambient space, e.g. $\he_\tI$ for fixed frame, $\be_\bI$ for moving frame.
    \item [Lowercase Greek characters] $\mu\in\{ 0,1,2,3 \}$ : \\
    The indices in the range of the basis of the submanifold, e.g. $\be_\bmu$, $\he_\tmu$.
    \item [Lowercase Latin characters] $i\in\{ 4,5,\dots,9 \}$ : \\
    The indices not belonging to the range of the basis of the submanifold, e.g. $\be_\bi$, $\he_\ti$.
\end{description}

With the above index notation, the standard metric of the ambient Minkowski space $\mathbb{R}^{1,9}$ is represented as $\eta_{\tI \tJ}$ (from the orthonormal fixed frame viewpoint), where $\tI,\tJ\in\{0,1,\dots,9\}$, and the metric of the pseudo-Riemannian manifold $\mathcal{M}^4$ is represented as $g_{\lmu\lnu}$, where  $\lmu,\lnu\in\{0,1,2,3\}$ and the underline font indices mean that the metric is viewed from the local coordinate system.

\subsection{Moving frame approach}

The Cartan moving frame approach \cite{Cartan:2001} describes the submanifold $\mathcal{M}$ and its tangent space by using the system $(\bx{q}(x), \{\be_\bI(x) \})$, where $\bx{q}(x)$ is the position vector $\vec{oq}$ of the point $q\in\mathcal{M}\subset\mathbb{R}^{1,9}$ (where $o\in\mathbb{R}^{1,9}$ is the origin point of the ambient space), and $\{\be_\bI(x)\}$ is the orthonormal moving frame on the position $q$ and $x$ is the local coordinate on the submanifold.
These vectors can be written in terms of the fixed frame, e.g. $\bx{q}(x) = q^{\tx{0}}(x)\hat{e}_{\tx{0}}+q^{\tx{1}}(x)\hat{e}_{\tx{1}} + \dots +q^{\tx{9}}(x)\hat{e}_\tx{9}$ which describes the local embedding map of $\bx{q}: \mathcal{M}\hookrightarrow\mathbb{R}^{1,9}$.
The following differential equation system describes the behaviour of the moving frame:
\begin{equation}\label{eq:moving_frame_equation}
    \left\{ \begin{array}{ccc}
        d\bx{q}(x) & = & \vartheta^\bI(x) \be_\bI(x) \\
        d\be_\bI(x) & = & \omega^{\;\;\bJ}_{\bI}(x) \be_{\bJ}(x)
    \end{array} \right.
\end{equation}
where $\vartheta^\bI(x) = \vartheta_{\lmu}^{\;\;\bI}(x) dx^{\lmu}$ is the vielbein 1-form field, and $\omega_\bI^{\;\; \bJ}(x) = \omega_{\lmu\bI}^{\;\;\;\;\bJ}(x) dx^{\lmu}$ is the connection 1-form on the submanifold $\mathcal{M}^4$, and here we use the Einstein summation convention for the repeated indices.
In fact, the index $\bI$ of the vielbein should only has the range $0\sim 3$, because the tangent space is spanned by $\be_{\bx{0}},\dots,\be_{\bx{3}}$.

If $g_{\lmu\lnu}$ is the metric of the local coordinate of submanifold, and $\bx{q}$ satisfies $ds^2 = g_{\lmu\lnu} dx^\lmu dx^\lnu = d\bx{q}(x)\cdot d\bx{q}(x) = \eta_{\bI \bJ} \vartheta_{\lmu}^{\;\;\bI}(x) \vartheta_{\lnu}^{\;\;\bJ}(x) dx^\lmu dx^\lnu$, it is called the isometric  immersion condition.
From the first order differential equation (\ref{eq:moving_frame_equation}), we know that the solution $\bx{q}(x)$ is deteremated by $\vartheta^\bI(x)$ and $\omega^{\;\;\bJ}_{\bI}(x)$ up to a initial condition.
And, the exact condition $dd\bx{q}=0$ shows $(d\vartheta^{\bI})\be_\bI - \vartheta^{\bK}\omega_{\bK}^{\;\;\bJ}\be_\bJ = 0$, it finally implies the immersion map $\bx{q}(x)$ can be deteremated just by the connection 1-form $\omega^{\;\;\bJ}_{\bI}(x)$ up to a initial condition.

Because the moving frame basis $\be_\bI$ are orthonormal $\be_\bI \cdot \be_\bJ=\eta_{\bI\bJ}$, so the connection 1-forms $\omega_{\bI\bJ}$ are antisymmetry:
\begin{align*}
    0 \;=\; d\left( \be_\bI \cdot \be_\bJ \right) \;=\; & d\be_\bI \cdot \be_\bJ + \be_\bI \cdot d \be_\bJ \\
    \;=\; & \omega^{\bK}_{\;\; \bI}\be_\bK \cdot \be_\bJ + \be_\bI \cdot \omega^{\bK}_{\;\; \bJ} \be_\bK \\
    \;=\; & \omega^{\bK}_{\;\; \bI}\eta_{\bK\bJ} + \omega^{\bK}_{\;\; \bJ}\eta_{\bI\bK} 
\end{align*}
i.e.
\begin{equation}
    \omega_{\bI\bJ} \;=\; -\omega_{\bJ\bI}
\end{equation}
According to the equation $d\be_\bI = \omega_{\bI}^{\;\; \bJ} \be_\bJ$, we can define the covariant derivative $D$
\begin{equation}
    D \be_\bI \;:=\; d \be_\bI - \omega_{\bI}^{\;\;\bJ} \be_\bJ \;=\; 0
\end{equation}
So, $\be_\bI$ are covariant constants under this covariant derivative.
We can define the curvature 2-form associated to this covariant derivative
\begin{equation}
    \Omega_{\bI}^{\;\;\bJ} \;:=\; D \omega_{\bI}^{\;\;\bJ} \;=\; d \omega_{\bI}^{\;\;\bJ} - \omega_{\bI}^{\;\;\bK} \wedge \omega_{\bK}^{\;\;\bJ}
\end{equation}
Due to $dd \be_\bI = 0$, it is easy to check that $\Omega_{\bI}^{\;\;\bJ}=0$.
We can seperate the indices into the tangent and normal subspace of the submanifold
\begin{equation}
    \left( \omega_{\bI}^{\;\; \bJ} \right) \;=\; \left( \begin{array}{cc}
        \omega_{\bmu}^{\;\; \bnu} & \omega_{\bmu}^{\;\; \bj} \\
        \omega_{\bi}^{\;\; \bnu} & \omega_{\bi}^{\;\; \bj}
    \end{array} \right)
\end{equation}
In order to easily distinguish them, we replace $\omega_{\bmu}^{\;\; \bj}$ with $H_{\bmu}^{\;\; \bj}$ , and call $H_{\bmu}^{\;\; \bj}$ the extrinsic curvature; and we replace $\omega_{\bi}^{\;\; \bj}$ with $A_{\bi}^{\;\; \bj}$, and call $A_{\bi}^{\;\; \bj}$ the connection of the normal bundle (or roughly call it the gauge field):
\begin{equation}
    \left( \omega_{\bI}^{\;\; \bJ} \right) \;=\; \left( \begin{array}{cc}
        \omega_{\bmu}^{\;\; \bnu} & H_{\bmu}^{\;\; \bj} \\
        H_{\bi}^{\;\; \bnu} & A_{\bi}^{\;\; \bj}
    \end{array} \right)
\end{equation}
The curvature 2-forms can be separately written in the tengent and normal subspace, these are the Gauss, Codazzi, Ricci equations:
\[
    \left\{ 
    \begin{array}{r l r}
        \textrm{Gauss eq. : } & \Omega_{\bmu}^{\;\; \bnu} \;=\; R_{\bmu}^{\;\; \bnu} - H_{\bmu}^{\;\; \bi} \wedge H_{\bi}^{\;\;\bnu} & \;=\; 0 \\
        \textrm{Codazzi eq. : } & \Omega_{\bmu}^{\;\; \bj} \;=\; d H_{\bmu}^{\;\; \bj} - \omega_{\bmu}^{\;\; \bsigma} \wedge H_{\bsigma}^{\;\; \bj}  - H_{\bmu}^{\;\; \bk} \wedge A_{\bk}^{\;\; \bj} & \;=\; 0\\
        \textrm{Ricci eq. : } & \Omega_{\bi}^{\;\; \bj} \;=\; F_{\bi}^{\;\; \bj} - H_{\bi}^{\;\; \bmu} \wedge H_{\bmu}^{\;\; \bj} & \;=\; 0
    \end{array}
    \right.
\]
where $R_{\bmu}^{\;\; \bnu} := d\omega_{\bmu}^{\;\;\bnu} - \omega_{\bmu}^{\;\;\bsigma} \wedge \omega_{\bsigma}^{\;\;\bnu}$ and $F_{\bi}^{\;\; \bj} := dA_{\bi}^{\;\; \bj} - A_{\bi}^{\;\; \bk} \wedge A_{\bk}^{\;\; \bj}$.

\subsection{Spin field approach}

To describe the spin group, we can use the Clifford algebra \cite{Shirokov:2017}.
The orthonormal basis $\{\he_\tI\}$ of the fixed frame in $\mathbb{R}^{1,9}$ will be seen as the generators of the Clifford algebra $Cl_{1,9}$, and they satisfies the rules $\he_\tI\he_\tJ = - \he_\tJ \he_\tI$ for $\tI\neq\tJ$ and $\he_\tI\he_\tJ+\he_\tJ\he_\tI = 2\eta_{\tI\tJ}$. 
The multivector set $\{ 1,\, \he_\tx{0},\, \he_\tx{1},\, \dots, \he_\tx{9},\, (\he_\tx{0}\he_\tx{1}),\, (\he_\tx{0}\he_\tx{2}), \dots, (\he_\tx{8}\he_\tx{9}),\; \dots, (\he_\tx{0}\he_\tx{1}\he_\tx{2}\dots) \}$ forms the basis of the Clifford algebra.
In the Clifford algebra $Cl_{1,9}$, the element $1 \in Cl_{1,9}^{0}$ is referred to as the element of grade 0, and $\he_\tI \in Cl_{1,9}^{1}$ is referred to as the element of grade 1, and $(\he_\tI\he_\tJ) \in Cl_{1,9}^{2}$ is referred to as the element of grade 2, and so on.
The even grade Clifford multivectors have the rotation property, and it is easy to check that $e^{\frac{-\theta}{2}\he_\tx{1}\he_\tx{2}}\he_\tx{1} e^{\frac{\theta}{2}\,\he_\tx{1}\he_\tx{2}} = \he_\tx{1} \cos(\theta) + \he_\tx{2} \sin(\theta)$, i.e. $\he_\tx{1}\he_\tx{2}$ is the generator of the rotation of the $(\he_\tx{1}-\he_\tx{2})$ plane.
In this article, a spin field $\psi$ is a even grade Clifford multivector function
\begin{equation}
    \psi= f + f^{\tI\tJ} \he_\tI\he_\tJ + f^{\tI\tJ\tK\tL} \he_\tI\he_\tJ\he_\tK\he_\tL + \dots
\end{equation}
, and satisfies the normalized condition $\tilde{\psi}\psi=1$, and $\tilde{\psi}\he_\tI\psi\in Cl_{1,9}^1$,\footnotemark[4] where the symbol $\sim$ is the reversion operation $\widetilde{AB}=\tilde{B}\tilde{A}$, i.e. $\widetilde{\he_\tx{1}\he_\tx{2}} = \he_\tx{2} \he_\tx{1}$.
\footnotetext[4]{
    We need this condition $\tilde{\psi}\he_\tI\psi\in Cl_{1,9}^1$ in the high dimensional space.
    Because, there is an example that even grade element $T = \frac{1}{\sqrt{2}}\left( \he_\tx{1}\he_\tx{2} + \he_\tx{3}\he_\tx{4}\he_\tx{5}\he_\tx{6} \right)\in Cl_n^{(even)}$ satisfies $\tilde{T}T=1$, but $\tilde{T}\he_\tx{1} T = \he_\tx{2}\he_\tx{3}\he_\tx{4}\he_\tx{5}\he_\tx{6}\neq Cl_n^1$.
    One can see the detail in \cite{Shirokov:2017}.
}

Due to the local rotation effect of the spin field, the deformed basis can be seen as a new moving frame basis
\begin{equation}
    \be_\bI \;=\; \tilde{\psi}\he_\tI\psi
\end{equation}
Note that at the moment we cannot guarantee that this $\be_\bI$ is the standard Cartan moving frame basis, because the differential equation can be more complex.
(In general, this kind of the moving frame basis can be $\be_\bJ = \sum_{\tI} a_{\bJ\tI}(x) \tilde{\psi}\he_\tI\psi$. 
However we can think that $\sum_{\tI} a_{\bJ\tI}(x) \he_\tI$ is another ortation effect and it can be absorbed into spin field.
Thus, in general, a moving frame basis can be written in the simple form $\be_\bI = \tilde{\psi}\he_\tI\psi$.)
Latter, we may abuse the text font indices and bold font indices, because their ranges are the same, and we will interpret the moving frame $\be_\bI$ as the deformed frame of $\he_\tI$ by a spin field $\psi$.
Physically, if one simply thinks the spin field is the spinor matter field, or the fermion.
Then, the basis deformation $\be_\bI = \tilde{\psi}\he_\tI\psi$ can be naively seen as an intuitive geometric interpretation that ``matter tells spacetime how to curve, and curved spacetime tells matter how to move''.

We will show that if a spin field $\psi$ satisfies the following Killing spin field equation
\begin{equation}\label{eq:killing_spin_field_equation}
    \partial_\lmu \psi \;=\; \frac{1}{4}\left( \sum_{\tnu=0}^3 \sum_{\tsigma=0}^3 \he_\tnu \he_\tsigma \omega^{\;\;\tnu\tsigma}_{\lmu} \right) \psi + \frac{1}{2}\left( \sum_{\tnu=0}^3 \sum_{\ti=4}^9 \he_\tnu \he_\ti H^{\;\;\tnu\ti}_{\lmu} \right) \psi + \frac{1}{4}\left( \sum_{\ti=4}^9 \sum_{\tj=4}^9 \he_\ti\he_\tj A^{\;\;\ti\tj}_{\lmu} \right) \psi
\end{equation}
, then $\tilde{\psi}\he_\tI\psi$ indeed has the same properties of the Cartan moving frame.
It is natural to request the antisymmetric conditions $\omega^{\;\;\tnu\tsigma}_{\lmu} = - \omega^{\;\;\tsigma\tnu}_{\lmu}$, $H^{\;\;\tnu\ti}_{\lmu} = - H^{\;\;\ti\tnu}_{\lmu}$, $A^{\;\;\ti\tj}_{\lmu} = -A^{\;\;\tj\ti}_{\lmu}$.
Late, we will use the Einstein summation convention and ignore the summation symbol.

Now, if a spin field $\psi$ satisfies the Killing spin field equation and the rotated basis is referred to as $\be_\bI = \tilde{\psi}\he_\tI\psi$, it is easy to check the differential relation has the same form as the Cartan moving frame
\begin{align*}
    d\be_\bmu \;=\; & \left( (\widetilde{\partial_\lsigma \psi})\he_\tmu \psi \right)dx^\lsigma + \left( \tilde{\psi}\he_\tmu (\partial_\lsigma \psi) \right) dx^\lsigma \\
    \;=\; & \frac{1}{4}\left( \tilde{\psi}\left( \omega^{\;\;\talpha\tbeta}_{\lsigma} \he_\tbeta \he_\talpha \right) + 2 \tilde{\psi}\left( H^{\;\;\tbeta\ti}_{\lsigma} \he_\ti\he_\tbeta \right) + \tilde{\psi} \left( A^{\;\;\ti\tj}_{\lsigma} \he_\tj \he_\ti \right) \right) \he_\tmu \psi dx^\lsigma \\
    & + \tilde{\psi} \he_\tmu \frac{1}{4}\Big( \left( \he_\talpha \he_\tbeta \omega^{\;\;\talpha\tbeta}_{\lsigma} \right) \psi + 2 \left( \he_\tbeta \he_\ti H^{\;\;\tbeta\ti}_{\lsigma} \right) \psi + \left( \he_\ti \he_\tj A^{\;\;\ti\tj}_{\lsigma} \right) \psi \Big) dx^\lsigma \\
    \;=\; & \omega^{\;\;\;\;\tbeta}_{\lsigma\tmu} \left( \tilde{\psi} \he_\tbeta \psi \right) dx^\lsigma + H^{\;\;\;\;\ti}_{\lsigma\tmu} \left( \tilde{\psi}\he_\ti \psi \right) dx^\lsigma \\
    \;=\; & \omega^{\;\;\;\;\bbeta}_{\lsigma\bmu} \be_\bbeta dx^\lsigma + H^{\;\;\;\;\bi}_{\lsigma\bmu} \be_\bi dx^\lsigma
\end{align*}
In the third equality above, we use the relation $\he_\tbeta \he_\talpha \he_\tmu = - \he_\tmu \he_\talpha \he_\tbeta$ for $\talpha\neq\tbeta\neq\tmu$, so the nonzero contributions come from the terms with contraction $\he_\tmu\he_\tnu=\eta_{\tmu\tnu}$ for $\tmu=\tnu$.
In the last equality above, we replace the index $\tmu$ by $\bmu$, it is allowed because the ranges are the same for $\tmu$ and $\bmu$.
Thus, by comparing with the equation of the moving frame, we see that the coefficient $\omega^{\;\;\;\;\bbeta}_{\lsigma\bmu}dx^\lsigma$ can be interpreted as the connection 1-form of the moving frame.
With the same reason, we get the similar relation for the normal basis vector
\begin{align*}
    d\be_\bi \;=\; & \left( (\widetilde{\partial_\lsigma \psi})\he_\ti \psi \right)dx^\lsigma + \left( \tilde{\psi}\he_\ti (\partial_\lsigma \psi) \right) dx^\lsigma \\
    \;=\; & \frac{1}{4}\left( \tilde{\psi}\left( \omega^{\;\;\talpha\tbeta}_{\lsigma} \he_\tbeta \he_\talpha \right) + 2\tilde{\psi}\left( H^{\;\;\tbeta\tj}_{\lsigma} \he_\tj \he_\tbeta \right) + \tilde{\psi} \left( A^{\;\;\tj\tk}_{\sigma} \he_\tk \he_\tj \right) \right) \he_\ti \psi dx^\lsigma \\
    & + \tilde{\psi} \he_\ti \frac{1}{4}\left( \left( \he_\talpha \he_\tbeta \omega^{\;\;\talpha\tbeta}_{\lsigma} \right) \psi + 2\left( \he_\tbeta \he_\tj H^{\;\;\tbeta\tj}_{\lsigma} \right) \psi + \left( \he_\tj \he_\tk A^{\;\;\tj\tk}_{\lsigma} \right) \psi \right) dx^\lsigma \\
    \;=\; & - H_{\lsigma\;\;\ti}^{\;\;\tbeta}\tilde{\psi}\he_\tbeta\psi dx^\lsigma + A_{\lsigma\ti}^{\;\;\;\;\tj} \tilde{\psi}\he_\tj\psi dx^\lsigma \\
    \;=\; & - H_{\lsigma\;\;\bi}^{\;\;\bbeta} \be_\bbeta dx^\lsigma + A_{\lsigma\bi}^{\;\;\;\;\bj} \be_\bj dx^\lsigma
\end{align*}
Thus, if $\psi$ satisfies the Killing spin field equation, then $\tilde{\psi}\he_\tI\psi$ has the same properties of the Cartan moving frame, and the coefficients have the direct interpretation.

And the integrability of the spin field $\psi$, i.e. $\partial_\lnu \partial_\lmu \psi=\partial_\lmu \partial_\lnu \psi$, will give the Gauss, Codazzi, Ricci equations.
\begin{small}
    \begin{align*}
        \partial_\lnu \partial_\lmu \psi \;=\; & \partial_\lnu \left[ \frac{1}{4}\left( \he_\talpha \he_\tbeta \omega^{\;\;\talpha\tbeta}_{\lmu} + 2 \he_\tbeta \he_\ti H_{\lmu}^{\;\;\tbeta\ti} + \he_\ti\he_\tj A^{\;\;\ti\tj}_{\lmu} \right) \psi \right] \\
        \;=\; & \frac{1}{4}\left( \he_\talpha \he_\tbeta \partial_\lnu \omega^{\;\;\talpha\tbeta}_{\lmu} + 2 \he_\tbeta \he_\ti \partial_\lnu H^{\;\;\tbeta\ti}_{\lmu} + \he_\ti\he_\tj \partial_\lnu A^{\;\;\ti\tj}_{\lmu} \right) \psi + \frac{1}{4}\left( \he_\talpha \he_\tbeta \omega^{\;\;\talpha\tbeta}_{\lmu} + 2 \he_\tbeta \he_\ti H^{\;\;\tbeta\ti}_{\lmu} + \he_\ti\he_\tj A^{\;\;\ti\tj}_{\lmu} \right) \partial_\lnu \psi \\
        \;=\; & \frac{1}{4}\left( \he_\talpha \he_\tbeta \partial_\lnu \omega^{\;\;\talpha\tbeta}_{\lmu} + \frac{1}{4} \he_\talpha \he_\tbeta \omega^{\;\;\talpha\tbeta}_{\lmu} \he_\tsigma \he_\tlambda \omega^{\;\;\tsigma\tlambda}_{\lnu} \right)\psi \\
        & + \frac{1}{4}\left( \he_\tbeta \he_\ti H^{\;\;\tbeta\ti}_{\lmu} \he_\tsigma \he_\tj H^{\;\;\tsigma\tj}_{\lnu} \right) \psi \\
        & + \frac{1}{4}\left( \he_\ti\he_\tj \partial_\lnu A_{\lmu}^{\;\;\ti\tj} + \frac{1}{4} \he_\ti\he_\tj A^{\;\;\ti\tj}_{\lmu} \he_\tk\he_\tl A^{\;\;\tk\tl}_{\lnu}  \right) \psi \\
        & + \frac{1}{8}\left( 4 \he_\tbeta \he_\ti \partial_\lnu H^{\;\;\tbeta\ti}_{\lmu} + \he_\ti\he_\tj A^{\;\;\ti\tj}_{\lmu} \he_\tbeta \he_\tk H^{\;\;\tbeta\tk}_{\lnu} + \he_\tbeta \he_\tk H^{\;\;\tbeta\tk}_{\lmu} \he_\ti\he_\tj A^{\;\;\ti\tj}_{\lnu} \right. \\
        & \qquad\qquad \left. + \he_\tsigma \he_\ti H^{\;\;\tsigma\ti}_{\lmu} \he_\talpha \he_\tbeta \omega^{\;\;\talpha\tbeta}_{\lnu} + \he_\talpha \he_\tbeta \omega^{\;\;\talpha\tbeta}_{\lmu} \he_\tsigma \he_\ti H^{\;\;\tsigma\ti}_{\lnu} \right)\psi \\
        & + \frac{1}{16} \left( \he_\talpha \he_\tbeta \omega^{\;\;\talpha\tbeta}_{\lmu} \he_\ti\he_\tj A^{\;\;\ti\tj}_{\lnu} + \he_\ti\he_\tj A^{\;\;\ti\tj}_{\lmu} \he_\talpha \he_\tbeta \omega^{\;\;\talpha\tbeta}_{\lnu} \right) \psi \\
        \;=\; & \frac{1}{4}\he_\talpha \he_\tbeta \left( \partial_\lnu \omega^{\;\;\talpha\tbeta}_{\lmu} + \omega^{\;\;\talpha}_{\lmu\;\;\tsigma} \omega^{\;\;\tsigma \tbeta}_{\lnu} + H^{\;\;\talpha\ti}_{\lmu} H^{\;\;\;\;\tbeta}_{\lnu\ti} \right) \psi \\
        & + \frac{1}{4} \he_\ti\he_\tj \left( \partial_\lnu A^{\;\;\ti\tj}_{\lmu} + A^{\;\;\ti}_{\lmu\;\;\tk} A^{\;\;\tk \tj}_{\lnu} + H^{\;\;\ti\tsigma}_{\lmu} H^{\;\;\;\;\tj}_{\lnu\tsigma}  \right) \psi \\
        & + \frac{1}{4} \he_\tbeta\he_\ti \left( 2 \partial_\lnu  H^{\;\;\tbeta\ti}_{\lmu} + A_{\lmu\;\;\tk}^{\;\;\ti} H^{\;\;\tbeta\tk}_{\lnu} - A_{\lnu\;\;\tk}^{\;\;\ti} H^{\;\;\tbeta\tk}_{\lmu} + \omega^{\;\;\tbeta}_{\lmu\;\;\tsigma} H^{\;\;\tsigma\ti}_{\lnu} - \omega^{\;\;\tbeta}_{\lnu\;\;\tsigma} H^{\;\;\tsigma\ti}_{\lmu} \right) \psi \\
        & + \he_\bullet \he_\bullet \he_\bullet \he_\bullet\left( \dots \right) \psi \; + \; 1 (\dots) \psi
    \end{align*}
\end{small}
where the last term $\he_\bullet \he_\bullet \he_\bullet \he_\bullet\left( \dots \right) \psi$ contains the rest of the multivectors of grade 4 which consist of the product of two Clifford algebras of grade 2, and all their vectors $\he$ are different, so the two Clifford algebras of grade 2 are commute.
Thus, for $\partial_\lmu\partial_\lnu\psi$ and $\partial_\lnu\partial_\lmu\psi$, the contents of the term $\he_\bullet \he_\bullet \he_\bullet \he_\bullet\left( \dots \right) \psi$ (, and also the term $1 (\dots) \psi$) will be the same.
Now, it is easy to see that the integrability $\partial_\lmu\partial_\lnu\psi = \partial_\lnu\partial_\lmu\psi$ gives the Gauss, Codazzi, Ricci equations.

\section{The solutions of the Killing spin field equation}

In this section, we discuss the solutions of the Killing spin field equation (\ref{eq:killing_spin_field_equation}).
Here, we assume that a spin field $\psi$ is written as
\begin{equation}
    \psi \;=\; f(x) + f^{\tI \tJ}(x) \he_\tI\he_\tJ + f^{\tI \tJ \tK \tL}(x) \he_\tI\he_\tJ\he_\tK\he_\tL + \dots
\end{equation}
, where $f(x), f^{\tI \tJ}(x), f^{\tI \tJ \tK \tL}(x),\dots$ are the functions of the submanifold coordinate $x^\lmu$.
First we notice that if $\psi$ is a solution of the Killing spin field equation (\ref{eq:killing_spin_field_equation}), then $(\he_\tx{4}\he_\tx{5}\dots\he_\tx{9})\psi$ is also a solution, but some coefficients need to be changed $H^{\tmu\tj}\to - H^{\tmu\tj}$.
For the same reason, in some simple cases, multiplying two solutions of the Killing spin field equation can generate new solutions, and we will discuss this later.
In the following subsections, we will give some nontrivial solutions of the Killing spin field equation.
We conjecture that all spin field representations of isometric immersion, i.e. all solutions of the Killing spin field equation, can be locally written as the products of the solutions we find.

It is easy to see that the rotations in the tangent space of the submanifold (, or the rotations in the normal space of the submanifold ) are the trivial solutions of the Killing spin field equation.
For example, $\psi = e^{\frac{1}{2}\theta\he_\tx{4}\he_\tx{5}}$, its partial derivative is $\partial_\lalpha\psi = \frac{1}{2}\left( \partial_\lalpha\theta \he_\tx{4}\he_\tx{5} \right) \psi$.
With the Killing spin field equation, it is easy to read out $A_\lalpha^{\;\;\tx{4}\tx{5}}=\partial_\lalpha\theta$, and its curvature tensor is zero $F_\ti^{\;\;\tj}=dA_\ti^{\;\;\tj}+A_\ti^{\;\;\tk}\wedge A_\tk^{\;\;\tj} = 0$, i.e. the associated submanifold is trivial.

\subsection{The type of $\psi = f + f^{\tmu \tx{5}}\he_\tmu\he_\tx{5}$}

Here, we show a spin field with the form 
\begin{equation}\label{eq:type_a_spin_field}
    \psi \;=\; f + f^{\tmu \tx{5}}\he_\tmu\he_\tx{5}
\end{equation}
and the normalized condition $\tilde{\psi}\psi=1$ will be a solution of the Killing spin field equation.
In order to take some explicit calculations, in equation (\ref{eq:type_a_spin_field}) we specify $\he_\tx{5}$ as the normal direction of the submanifold, however it can be replaced by any one of the normal directions $\he_\ti$.
Note that the repeated symbolic index $\tmu$ obeys the Einstein summation convention, and there is only one normal direction of the submanifold in the spin field.

We substitute the spin field (\ref{eq:type_a_spin_field}) into the Killing spin field equation (\ref{eq:killing_spin_field_equation}), and expand it in terms of the different multivector basis.
Then we get the following equation set:
\begin{equation}\label{eq:set_of_killing_spin_field_equations}
    \begin{array}{r c r c l}
        1 & : \qquad & \partial_\lalpha f & = & - \frac{1}{2} H_{\lalpha\tmu\tx{5}} f^{\tmu \tx{5}} \\
        \he_\tmu\he_\tx{5} & : \qquad & \partial_\alpha f^{\tmu \tx{5}} & = & \frac{1}{2}\omega^{\;\;\tmu}_{\lalpha\;\;\tnu} f^{\tnu \tx{5}} + \frac{1}{2}H^{\;\;\tmu\tx{5}}_{\lalpha}f \\
        \he_\tmu\he_\tnu & : \qquad & 0 = \partial_\lalpha f^{\tmu\tnu} & = & \frac{1}{2}\omega^{\;\;\tmu\tnu}_{\lalpha} f - \frac{1}{2}H^{\;\;\tmu}_{\lalpha\;\;\tx{5}}f^{\tnu \tx{5}} + \frac{1}{2} H^{\;\;\tnu}_{\lalpha\;\;\tx{5}}f^{\tmu \tx{5}} \\
        \he_\tx{5}\he_\tx{i} & : \qquad & 0 = \partial_\lalpha f^{\tx{5}\ti} & = & \frac{1}{2} H^{\;\;\;\;\ti}_{\lalpha\tmu}f^{\tmu \tx{5}} + \frac{1}{2} A^{\;\;\tx{5}\ti}_{\lalpha}f \\
        & & \vdots
    \end{array}
\end{equation}
, where the upper (or lower) indices of the text fonts are lowered (or raised) by the metric $\eta_{\tI\tJ}$.
The third equation of the equation set (\ref{eq:set_of_killing_spin_field_equations}) says 
\begin{equation}\label{eq:connection_in_terms_of_extrinsic_curvature}
    \omega^{\;\;\tmu\tnu}_{\lalpha} = \frac{1}{f}\left( H^{\;\;\tmu}_{\lalpha\;\;\tx{5}}f^{\tnu \tx{5}} - H^{\;\;\tnu}_{\lalpha\;\;\tx{5}}f^{\tmu \tx{5}} \right)
\end{equation}
Substituting this result into the second equation of the equation set (\ref{eq:set_of_killing_spin_field_equations}), we get $\partial_\lalpha f^{\tmu \tx{5}} = \frac{1}{2}\frac{1}{f}\left( H^{\;\;\tmu\tx{5}}_{\lalpha}f_{\tnu\tx{5}} - H_{\lalpha\tnu\tx{5}}f^{\tmu\tx{5}} \right) f^{\tnu \tx{5}} + \frac{1}{2}H^{\;\;\tmu\tx{5}}_{\lalpha} f$.
Finally, using the first equation of the equation set (\ref{eq:set_of_killing_spin_field_equations}) and the normalized condition $(f)^2 + f^{\tmu \tx{5}} f^{\;\;\tx{5}}_\tmu =1$, we get
\[
    \partial_\lalpha f^{\tmu \tx{5}} \;=\; \frac{1}{2}\frac{1}{f}H^{\;\;\tmu\tx{5}}_{\lalpha} + \frac{1}{f}\left( \partial_\lalpha f \right) f^{\tmu \tx{5}}
\]
Thus, the extrinsic curvature is
\begin{equation}
    H^{\;\;\tmu\tx{5}}_{\lalpha} \;=\; 2 (f)^2 \partial_\lalpha\left( \frac{f^{\tmu \tx{5}}}{f} \right)
\end{equation}
This result and eq.(\ref{eq:connection_in_terms_of_extrinsic_curvature}) show that if there are two terms $f^{\tmu\tx{5}}$ non-zero, then $\omega_{\lalpha}^{\;\;\tmu\tnu} \neq 0$ and the curvature should be non-zero.

Another equation that is not shown in equation set (\ref{eq:set_of_killing_spin_field_equations}) is $0=\partial_\lalpha f^{\tmu\ti} = \frac{1}{2}H_\lalpha^{\;\;\tmu\ti}f - \frac{1}{2}A_{\lalpha}^{\;\;\tx{5}\ti}f^{\tmu\tx{5}}$ for $\ti\neq\tx{5}$, and it implies $H^{\;\;\tmu\ti}_{\lalpha} = \frac{1}{f}A_{\lalpha}^{\;\;\tx{5}\ti}f^{\tmu\tx{5}}$ for $\ti\neq \tx{5}$.
With the result, the fourth equation of the equation set (\ref{eq:set_of_killing_spin_field_equations}) implies $A^{\;\;\tx{5}\ti}_{\lalpha}=0$ and $H^{\;\;\tmu\ti}_{\lalpha}=0$ for $\ti\neq \tx{5}$.
And, one can check that this solution will automately satisfy the equation $0=\partial_\lalpha f^{\tmu\tnu\tgamma \tx{5}} = \frac{1}{4}\omega^{\;\;\left[\tmu\tnu \right.}_{\lalpha}f^{\left. \tgamma \right] \tx{5}}$.
Finally, it has been shown that the field (\ref{eq:type_a_spin_field}) is the solution of the Killing spin field equation.

\subsubsection{Example}

We consider a simple example, the spin field is
\begin{equation}
    \psi \;=\; \frac{1}{\sqrt{1+r^2}} - \frac{x_\lx{1}}{\sqrt{1+r^2}}\he_\tx{1}\he_\tx{5} - \frac{x_\lx{2}}{\sqrt{1+r^2}}\he_\tx{2}\he_\tx{5} - \frac{x_\lx{3}}{\sqrt{1+r^2}}\he_\tx{3}\he_\tx{5}
\end{equation}
where $r^2=x_\lx{1}^2+x_\lx{2}^2+x_\lx{3}^2$.
It is easy to check $\tilde{\psi}\psi=1$.

Using the formulae in the previous subsection, we obtain the associated extrinsic curvature: 
\[
    \begin{array}{c c l}
        H^{\;\;\bx{0}\bx{5}}_{\lx{0}}\;=\;0 \; , \quad H^{\;\;\bmu\bx{5}}_{\lmu} \;=\; \frac{-2}{(1+r^2)} & \quad & \textrm{for $\bmu=\lmu\neq 0$} \\
        H^{\;\;\bmu\bx{5}}_{\lnu}\;=\;0 & & \textrm{for $\bmu\neq \lnu$}
    \end{array}
\]
, and the connection: 
\[
    \begin{array}{l l l}
        \omega^{\;\;\bmu\bnu}_{\lx{0}}\;=\;0 \;,\; & \omega^{\;\;\bx{0}\bmu}_{\lnu}=0 & \\
        \omega^{\;\;\bx{1}\bx{2}}_{\lx{1}}\;=\;\frac{2 x_\lx{2}}{(1+r^2)} \;,\; & \omega^{\;\;\bx{1}\bx{3}}_{\lx{1}}\;=\;\frac{2 x_\lx{3}}{(1+r^2)} \;,\; & \omega^{\;\;\bx{2}\bx{3}}_{\lx{1}}\;=\;0 \\
        \omega^{\;\;\bx{1}\bx{2}}_{\lx{2}}\;=\;\frac{-2 x_\lx{1}}{(1+r^2)} \;,\; & \omega^{\;\;\bx{1}\bx{3}}_{\lx{2}}\;=\;0 \;,\; & \omega^{\;\;\bx{2}\bx{3}}_{\lx{2}}\;=\;\frac{2 x_\lx{3}}{(1+r^2)} \\
        \omega^{\;\;\bx{1}\bx{2}}_{\lx{3}}\;=\;0 \;,\; & \omega^{\;\;\bx{1}\bx{3}}_{\lx{3}}\;=\;\frac{-2 x_\lx{1}}{(1+r^2)} \;,\; & \omega^{\;\;\bx{2}\bx{3}}_{\lx{3}}\;=\;\frac{-2 x_\lx{2}}{(1+r^2)}
    \end{array}
\]
Thus, the curvature tensors ($R_{\bmu}^{\;\;\bnu}=d\omega_{\bmu}^{\;\;\bnu} - \omega_{\bmu}^{\;\;\bsigma}\wedge\omega_{\bsigma}^{\;\;\bnu}$) are
\begin{align}
    & R^{\;\;\;\;\bx{1}\bx{2}}_{\lx{1}\lx{2}} \;=\; R^{\;\;\;\;\bx{1}\bx{3}}_{\lx{1}\lx{3}} \;=\; R^{\;\;\;\;\bx{2}\bx{3}}_{\lx{2}\lx{3}} \;=\; - \frac{4}{(1+r^2)^2} \label{eq:curvature_tensors}\\
    & \textrm{other $R^{\;\;\;\;\bmu\bnu}_{\lalpha\lbeta}=0$} \nonumber
\end{align}
In the above curvature tensor, the uper indices describe the moving frame $\be_\balpha$ and the lower indices describe the local coordinate $x_\lmu$ of the manifold $\mathcal{M}^4$.

To find the local (pseudo-)Riemannian metric $g_{\lmu\lnu}$, we can try to find the solution of the moving frame equation $d\bx{q}=\vartheta^\bI\be_\bI$.
The basis of the moving frame are $\be_\bI=\tilde{\psi}\he_\tI\psi$ :
\begin{align*}
    \be_\bx{0} \;=\; & \he_0 \\
    \be_\bx{1} \;=\; & \frac{1-x_\lx{1}^2+x_\lx{2}^2+x_\lx{3}^2}{1+r^2}\he_\tx{1} - \frac{2x_\lx{1}x_\lx{2}}{1+r^2}\he_\tx{2} - \frac{2x_\lx{1}x_\lx{3}}{1+r^2}\he_\tx{3} - \frac{2x_\lx{1}}{1+r^2}\he_\tx{5} \\
    \be_\bx{2} \;=\; & - \frac{2x_\lx{1}x_\lx{2}}{1+r^2}\he_\tx{1} + \frac{1+x_\lx{1}^2-x_\lx{2}^2+x_\lx{3}^2}{1+r^2}\he_\tx{2} - \frac{2x_\lx{2}x_\lx{3}}{1+r^2}\he_\tx{3} - \frac{2x_\lx{2}}{1+r^2}\he_\tx{5} \\
    \be_\bx{3} \;=\; & - \frac{2x_\lx{1}x_\lx{3}}{1+r^2}\he_\tx{1} + \frac{2x_\lx{2}x_\lx{3}}{1+r^2}\he_\tx{2} + \frac{1+x_\lx{1}^2+x_\lx{2}^2-x_\lx{3}^2}{1+r^2}\he_\tx{3} - \frac{2x_\lx{3}}{1+r^2}\he_\tx{5} \\
    \be_\bx{5} \;=\; & \frac{2x_\lx{1}}{1+r^2}\he_\tx{1} + \frac{2x_\lx{2}}{1+r^2}\he_\tx{2} + \frac{2x_\lx{3}}{1+r^2}\he_\tx{3} + \left( \frac{2}{1+r^2} - 1 \right)\he_\tx{5}
\end{align*}
For this simple case, it is not hard to guess the isometric immersion map is $\bx{q}=x_\lx{0}\he_\tx{0} + \frac{x_\lx{1}}{1+r^2}\he_\tx{1}+\frac{x_\lx{2}}{1+r^2}\he_\tx{2}+\frac{x_\lx{3}}{1+r^2}\he_\tx{3}+\frac{1}{1+r^2}\he_\tx{5}$, and the vielbein $\vartheta^{\;\;\bmu}_{\lnu}={\rm diag}(1,\frac{1}{1+r^2},\frac{1}{1+r^2},\frac{1}{1+r^2})$, $\vartheta^{\;\;\tx{5}}_\lnu=0$; they satisfy the relation $d\bx{q}=\vartheta^\bI\be_\bI$.
The isometric  immersion condition $g_{\lmu\lnu} = \partial_\lmu \bx{q} \cdot \partial_\lnu \bx{q}$ gives $g_{\lmu\lnu} = {\rm diag}(-1,\frac{1}{(1+r^2)^2},\frac{1}{(1+r^2)^2},\frac{1}{(1+r^2)^2})$.
Using the metric, we can calculate the Riemannian curvature tensor of the tangent bundle and check it is consistent with the result (\ref{eq:curvature_tensors}).

\subsubsection{Generating new solutions}\label{sec:generating_new_solutions1}

Assume $\psi_1$, $\psi_2$ are two spin fields, which satisfy the Killing spin field equation $\partial \psi_1 = K_1\psi_1$, $\partial \psi_2 = K_2\psi_2$, where $K_1$, $K_2$ only contain the Clifford algebras of grade 2 as in eq.(\ref{eq:killing_spin_field_equation}), then
\begin{align*}
    \partial \left( \psi_1\psi_2 \right) \;=\; & K_1\psi_1\psi_2 + \psi_1 K_2\psi_2 \\
    \;=\; & K_1\psi_1\psi_2 + \psi_1K_2\tilde{\psi}_1 \psi_1\psi_2 \\
    \;=\; & K_1\psi_1\psi_2 + K_3 \psi_1\psi_2
\end{align*}
It is easy to check that if $\psi_1$ only contains the Clifford algebras of grade 0 and grade 2, then $K_3=\psi_1K_2\tilde{\psi}_1$ is also the function of the Clifford algebra of grade 2.
Thus, $(\psi_1\psi_2)$ is also a solution of the Killing spin field equation, and the associated connections are the combinations of $K_1$ and $K_3$.

Here we give the explicit result for $\psi_1=f + \he_\tmu\he_\tx{5}f^{\tmu \tx{5}}$ and
\begin{align*}
    \partial_\lalpha\psi_1 \;=\; & \left( \frac{1}{4}\he_\tmu\he_\tnu \omega^{{}^{(1)}\;\;\tmu\tnu}_{\;\;\;\;\lalpha} + \frac{1}{2}\he_\tmu\he_\tx{5} H^{^{(1)}\;\;\tmu\tx{5}}_{\;\;\;\;\lalpha} \right)\psi_1 \\
    \partial_\lalpha\psi_2 \;=\; & \left( \frac{1}{4}\he_\tmu\he_\tnu\omega^{{}^{(2)}\;\;\tmu\tnu}_{\;\;\;\;\lalpha} + \frac{1}{2}\he_\tmu\he_\ti H^{{}^{(2)}\;\;\tmu\ti}_{\;\;\;\;\lalpha} + \frac{1}{4}\he_\ti\he_\tj A^{{}^{(2)}\;\; \ti\tj}_{\;\;\;\;\lalpha} \right)\psi_2 \\
\end{align*}
One can calculate the following equation, and we leave the detailed calculation to the appendix A.
\begin{align*}
    \partial_\lalpha\left( \psi_1\psi_2 \right) \;=\; & \left( \frac{1}{4}\he_\tmu\he_\tnu\omega^{{}^{(1)}\;\;\tmu\tnu}_{\;\;\;\;\lalpha} + \frac{1}{2}\he_\tmu\he_\tx{5} H^{{}^{(1)}\;\;\tmu\tx{5}}_{\;\;\;\;\lalpha} \right)\psi_1\psi_2 \\
    & + \left( \psi_1\left( \frac{1}{4}\he_\tmu\he_\tnu\omega^{{}^{(2)}\;\;\tmu\tnu}_{\;\;\;\;\lalpha} + \frac{1}{2}\he_\tmu\he_\ti H^{{}^{(2)}\;\;\tmu\ti}_{\;\;\;\;\lalpha} + \frac{1}{4}\he_\ti\he_\tj A^{{}^{(2)}\;\; \ti\tj}_{\;\;\;\;\lalpha} \right) \tilde{\psi}_1 \right) \psi_1 \psi_2 \\
    \;=\; & \left( \frac{1}{4}\he_\tmu\he_\tnu {\omega'}^{\;\;\tmu\tnu}_\lalpha + \frac{1}{2}\he_\tmu\he_\ti {H'}^{\;\;\tmu\ti}_\lalpha + \frac{1}{4}\he_\ti\he_\tj {A'}^{\;\; \ti\tj}_\lalpha \right) \psi_1 \psi_2
\end{align*}
Then, one can find that for the spin field $(\psi_1\psi_2)$, the associated connections ${\omega'}^{\;\;\tmu\tnu}_\lalpha$ are
\begin{align*}
    {\omega'}^{\;\;\tmu\tnu}_{\lalpha} \;=\; & \omega^{{}^{(1)}\;\;\tmu\tnu}_{\;\;\;\;\lalpha} + \omega^{{}^{(2)}\;\;\tmu\tnu}_{\;\;\;\;\lalpha}\left( (f)^2 + \sum_{\tgamma\neq \tmu,\tnu} f_{\tgamma \tx{5}}f^{\tgamma \tx{5}} \right) + 2\omega^{{}^{(2)}\;\;\tnu}_{\;\;\;\;\lalpha\;\;\tgamma}(f^{\tmu \tx{5}}f^{\tgamma \tx{5}}) - 2\omega^{{}^{(2)}\;\;\tmu}_{\;\;\;\;\lalpha\;\;\tgamma}(f^{\tnu}_{\;\;\tx{5}}f^{\tgamma \tx{5}}) \\
    & + 2 H^{{}^{(2)}\;\;\tmu}_{\;\;\;\;\lalpha\;\;\tx{5}} \left( f\,f^{\tnu \tx{5}} \right) - 2 H^{{}^{(2)}\;\;\tnu}_{\;\;\;\;\lalpha\;\;\tx{5}} \left( f\,f^{\tmu \tx{5}} \right)
\end{align*}
, the associated extrinsic curvatures are
\begin{align*}
    {H'}^{\;\;\tmu\tx{5}}_{\lalpha} \;=\; & H^{{}^{(2)}\;\;\tmu\tx{5}}_{\;\;\;\;\lalpha}\left( (f)^2 - 2\sum_{\tgamma\neq\tmu} f_{\tgamma \tx{5}}f^{\tgamma \tx{5}} \right) - 2 \omega^{{}^{(2)}\;\;\tmu}_{\;\;\;\;\lalpha\;\;\tgamma}(f\,f^{\tgamma \tx{5}}) + 4 H^{{}^{(2)}\;\;\;\;}_{\;\;\;\;\lalpha\tgamma\tx{5}}(f^{\tgamma \tx{5}}f^{\tmu \tx{5}}) \\
    {H'}^{\;\;\tmu\ti}_{\lalpha} \;=\; & H^{{}^{(2)}\;\;\tmu\ti}_{\;\;\;\;\lalpha}\left( (f)^2 + \sum_{\tgamma\neq\tmu} f_{\tgamma \tx{5}}f^{\tgamma \tx{5}} \right) + 2 A^{{}^{(2)}\;\;\;\;\ti}_{\;\;\;\;\lalpha\tx{5}} \left( f\,f^{\tmu \tx{5}} \right) + 2 H^{{}^{(2)}\;\;\;\;\ti}_{\;\;\;\;\lalpha\tgamma}(f^{\tgamma}_{\;\;\tx{5}}f^{\tmu \tx{5}}) \quad,\quad \ti\neq\tx{5}
\end{align*}
, and the associated gauge field are
\begin{align*}
    {A'}_\lalpha^{\;\;\tx{5}\ti} \;=\; & A^{{}^{(2)}\;\;\tx{5}\ti}_{\;\;\;\;\lalpha} \left( (f)^2 - 2 f_{\tgamma \tx{5}}f^{\tgamma \tx{5}} \right) - 2 H^{{}^{(2)}\;\;\;\;\ti}_{\;\;\;\;\lalpha\tgamma}(f\,f^{\tgamma \tx{5}}) \\
    {A'}_\lalpha^{\;\;\ti\tj} \;=\; & A^{{}^{(2)}\;\;\ti\tj}_{\;\;\;\;\lalpha}  \quad,\quad \ti,\tj\neq\tx{5}
\end{align*}
where the repeated symbolic indices obey the Einstein summation convention, except the specifying summation range with the summation symbol.

From the relation of ${\omega'}^{\;\;\tmu\tnu}_{\lalpha}$, we find that it is linear dependent on $\omega^{{}^{(1)}\;\;\tmu\tnu}_{\;\;\;\;\lalpha}$, $\omega^{{}^{(2)}\;\;\tmu\tnu}_{\;\;\;\;\lalpha}$, and $H^{{}^{(2)}\;\;\tmu\tx{5}}_{\;\;\;\;\lalpha}$.
And the nonlinear effects come form spin field $f\,f^{\tmu \tx{5}}$ and their squares.

From the relation of ${A'}_\lalpha^{\;\;\tx{5}\ti}$, we find that the magnitude and frequency of the gauge field can be changed.
If we really interpret the field $A$ as the physical gauge field.
This formula should be useful for studying the light redshift in the externally applied dynamical space-time (or in the gravitational wave background).

{\bf Remark:} Here, we just generate some new submanifolds locally, and don't request any physical rules.
There are some possible further research directions for the applications in physics : 
\begin{enumerate}[label=(\arabic*)]
    \item We may try to find a least action principle for the spin fields that is equivalent to the Einstein-Hilbert action.
    \item Using the above equations, we may study the stable condition. For example, we can ask what kind of connection $\omega^{{}^{(2)}\;\;\tmu\tnu}_{\;\;\;\;\lalpha}$ and extrinsic curvature $H^{{}^{(2)}\;\;\tmu\ti}_{\;\;\;\;\lalpha}$ are stable under small spin field $(f\approx 1 \;\&\; f^{\tmu\tx{5}}\approx \epsilon)$ perturbation ?
\end{enumerate}

\subsection{The type of $\psi = f + f^{\tx{1} \tk}\he_\tx{1}\he_\tk$}

In this section we show a spin field with the form
\begin{equation}
    \psi \;=\; f + f^{\tx{1} \tk}\he_\tx{1}\he_\tk
\end{equation}
will be a solution of the Killing spin field equation.
In order to take some explicit calculations, we specify $\he_\tx{1}$ as the tangent direction of
the submanifold, however in general it can be replaced by any one of the tangent directions $\he_\tmu$.
Note that the repeated symbolic index $\tk$ index obeys the Einstein summation convention, and there is only one tangent direction of the submanifold in the spin field.

We substitute the spin field into the Killing spin field equation (\ref{eq:killing_spin_field_equation}), and expand
it in terms of the different multivector basis.
\begin{equation}
    \begin{array}{r c r c l}
        1 & : \qquad & \partial_\lalpha f & = & - \frac{1}{2} H_{\lalpha\tx{1}\tk} f^{\tx{1}\tk} \\
        \he_\tx{1}\he_\tk & : \qquad & \partial_\alpha f^{\tx{1} \tk} & = & \frac{1}{2}H^{\;\;\tx{1}\tk}_{\lalpha}f + \frac{1}{2} A^{\;\;\tk}_{\lalpha\;\;\tj} f^{\tx{1}\tj} \\
        \he_\tx{1}\he_\tmu & : \qquad & 0 = \partial_\lalpha f^{\tx{1}\tmu} & = & \frac{1}{2}\omega^{\;\;\tx{1}\tmu}_{\lalpha} f + \frac{1}{2}H^{\;\;\tmu}_{\lalpha\;\;\tk}f^{\tx{1}\tk} \\
        \he_\ti\he_\tj & : \qquad & 0 = \partial_\lalpha f^{\ti\tj} & = & - \frac{1}{2} H^{\;\;\;\;\ti}_{\lalpha\tx{1}}f^{\tx{1} \tj} + \frac{1}{2} H^{\;\;\;\;\tj}_{\lalpha\tx{1}}f^{\tx{1} \ti} + \frac{1}{2} A^{\;\;\ti\tj}_{\lalpha}f \\
        & & \vdots
    \end{array}
\end{equation}
With the similar calculation in the previous section, one can get
\begin{equation}
    A_{\lalpha}^{\;\;\ti\tj} \;=\; \frac{1}{f}\left( H_{\lalpha\tx{1}}^{\;\;\;\;\ti}f^{\tx{1}\tj} - H_{\lalpha\tx{1}}^{\;\;\;\;\tj}f^{\tx{1}\ti} \right)
\end{equation}
and
\begin{equation}
    H_{\lalpha}^{\;\;\tx{1}\tk} \;=\; 2 (f)^2\partial_\lalpha\left( \frac{f^{\tx{1}\tk}}{f} \right)
\end{equation}

In this case, the connections $\omega_\lalpha^{\;\;\tmu\tnu}=0$ are zero, this means that these submanifolds made entirely of $\psi = f + f^{\tx{1} \tk}\he_\tx{1}\he_\tk$ are flat, so they are diffeomorphic to $\mathbb{R}^{1,3}$.

In next section we will see the spin field $\psi = f + f^{\tx{1} \tk}\he_\tx{1}\he_\tk$ can actually deform the submanifold that the connections or extrinsic curvatures are not zeros $\omega_\lalpha^{\;\;\tmu\tnu}\neq 0$ or $H_\lalpha^{\;\;\tmu\ti}\neq 0$.
(If we consider the spin field is a fermion $\psi$ that satisfies $\partial_\lalpha\psi = \left( \frac{1}{2}\he_\tmu\he_\ti H^{\;\;\tmu\ti}_{\lalpha} + \frac{1}{4}\he_\ti\he_\tj A^{\;\;\ti\tj}_{\lalpha} \right)\psi$.
From GR viewpoint, the fermions should deform the geometry of the submanifold.
So, the result in the next section can be interpreted as matter deforming the space-time.)

\subsubsection{Generating new solutions}\label{sec:generating_new_solutions2}

Similar with section {\ref{sec:generating_new_solutions1}}, here let $\psi_1= f + \he_\tx{1} \he_\tk f^{\tx{1}\tk}$ and
\begin{align*}
    \partial_\lalpha\psi_1 \;=\; & \left( \frac{1}{2}\he_\tx{1}\he_\tk H^{{}^{(1)}\;\;\tx{1}\tk}_{\;\;\;\;\lalpha} + \frac{1}{4} \he_\ti \he_\tj A^{{}^{(1)\;\;\ti\tj}}_{\;\;\;\;\lalpha} \right) \psi_1 \\
    \partial_\lalpha\psi_2 \;=\; & \left( \frac{1}{4}\he_\tmu\he_\tnu\omega^{{}^{(2)}\;\;\tmu\tnu}_{\;\;\;\;\lalpha} + \frac{1}{2}\he_\tmu\he_\ti H^{{}^{(2)}\;\;\tmu\ti}_{\;\;\;\;\lalpha} + \frac{1}{4}\he_\ti\he_\tj A^{{}^{(2)}\;\; \ti\tj}_{\;\;\;\;\lalpha} \right)\psi_2 \\
\end{align*}
then
\begin{align*}
    \partial_\lalpha\left( \psi_1\psi_2 \right) \;=\; & \left( \frac{1}{2}\he_\tx{1}\he_\tk H^{{}^{(1)}\;\;\tx{1}\tk}_{\;\;\;\;\lalpha} + \frac{1}{4} \he_\ti \he_\tj A^{{}^{(1)\;\;\ti\tj}}_{\;\;\;\;\lalpha} \right) \psi_1\psi_2 \\
    & + \left( \psi_1\left( \frac{1}{4}\he_\tmu\he_\tnu\omega^{{}^{(2)}\;\;\tmu\tnu}_{\;\;\;\;\lalpha} + \frac{1}{2}\he_\tmu\he_\ti H^{{}^{(2)}\;\;\tmu\ti}_{\;\;\;\;\lalpha} + \frac{1}{4}\he_\ti\he_\tj A^{{}^{(2)}\;\; \ti\tj}_{\;\;\;\;\lalpha} \right) \tilde{\psi}_1 \right) \psi_1 \psi_2 \\
    \;=\; & \left( \frac{1}{4}\he_\tmu\he_\tnu {\omega'}^{\;\;\tmu\tnu}_\lalpha + \frac{1}{2}\he_\tmu\he_\ti {H'}^{\;\;\tmu\ti}_\lalpha + \frac{1}{4}\he_\ti\he_\tj {A'}^{\;\; \ti\tj}_\lalpha \right) \psi_1 \psi_2
\end{align*}
One can find that for the spin field $(\psi_1\psi_2)$ the associated connections ${\omega'}^{\;\;\tmu\tnu}_\lalpha$ are
\begin{align*}
    {\omega'}^{\;\;\tx{1}\tmu}_{\lalpha} \;=\; & \omega^{{}^{(2)}\;\;\tx{1}\tmu}_{\;\;\;\;\lalpha}(f)^2 - 2 \omega^{{}^{(2)}\;\;\;\;\tmu}_{\;\;\;\;\lalpha\tx{1}}\left( f^{\tx{1}}_{\;\;\tk}f^{\tx{1}\tk} \right) - 2 H^{{}^{(2)}\;\;\tmu}_{\;\;\;\;\lalpha\;\;\tk}(f\,f^{\tx{1}\tk}) \\
    {\omega'}^{\;\;\tmu\tnu}_{\lalpha} \;=\; & \omega^{{}^{(2)}\;\;\tmu\tnu}_{\;\;\;\;\lalpha} \qquad,\qquad\qquad \textrm{for $\tmu,\tnu\neq \tx{1}$}
\end{align*}
, the associated extrinsic curvatures are
\begin{align*}
    {H'}^{\;\;\tx{1}\ti}_{\lalpha} \;=\; & H^{{}^{(1)}\;\;\tx{1}\ti}_{\;\;\;\;\lalpha} + H^{{}^{(2)}\;\;\tx{1}\ti}_{\;\;\;\;\lalpha} \left( ( f )^2 - 2 \sum_{\tk\neq\ti} f_{\tx{1} \tk}f^{\tx{1} \tk} \right) - 2 A^{{}^{(2)}\;\; \ti}_{\;\;\;\;\lalpha\;\;\tk} \left( f\,f^{\tx{1} \tk} \right) + 4 H^{{}^{(2)}\;\;}_{\;\;\;\;\lalpha\tx{1}\tk}\left( f^{\tx{1} \tk}f^{\tx{1} \ti} \right) \\
    {H'}^{\;\;\tmu\ti}_{\lalpha} \;=\; & H^{{}^{(2)}\;\;\tmu\ti}_{\;\;\;\;\lalpha} \left( ( f )^2 + \sum_{\tk\neq\ti} f_{\tx{1} \tk}f^{\tx{1} \tk} \right) - 2 \omega^{{}^{(2)}\;\;\tmu}_{\;\;\;\;\lalpha\;\;\tx{1}} \left( f\,f^{\tx{1} \ti} \right) - 2 H^{{}^{(2)}\;\;\tmu}_{\;\;\;\;\lalpha\;\;\tk}\left( f_\tx{1}^{\;\; \ti}f^{\tx{1} \tk} \right) \quad,\quad\textrm{for $\tmu\neq\tx{1}$}
\end{align*}
, and the associated gauge field are
\begin{align*}
    {A'}^{\;\;\ti\tj}_\lalpha \;=\; & A^{{}^{(1)}\;\;\ti\tj}_{\;\;\;\;\lalpha} + A^{{}^{(2)}\;\;\ti\tj}_{\;\;\;\;\lalpha} \left( ( f )^2 + \sum_{\tk\neq \ti,\tj} f_{\tx{1} \tk}f^{\tx{1} \tk} \right) + 2 A^{{}^{(2)}\;\;\tj}_{\;\;\;\;\lalpha\;\;\tm}\left( f_\tx{1}^{\;\;\ti}f^{\tx{1}\tm} \right) - 2 A^{{}^{(2)}\;\;\ti}_{\;\;\;\;\lalpha\;\;\tm}\left( f_\tx{1}^{\;\;\tj}f^{\tx{1}\tm} \right) \\
    & \qquad + 2 H^{{}^{(2)}\;\;\;\;\ti}_{\;\;\;\;\lalpha\tx{1}} \left( f\,f^{\tx{1} \tj} \right) - 2 H^{{}^{(2)}\;\;\;\;\tj}_{\;\;\;\;\lalpha\tx{1}} \left( f\,f^{\tx{1} \ti} \right)
\end{align*}
where the repeated symbolic indices obey the Einstein summation convention, except the specifying summation range with the summation symbol.

\subsection{A conjecture for the general solutions}

\textbf{Conjecture:} Any spin field $\psi$ which satisfies the Killing spin field equation (\ref{eq:killing_spin_field_equation}) locally can be written as the product of the first type $(f_{{}_{(4)}}+f_{{}_{(4)}}^{\;\;\tmu\tx{4}}\he_\tmu\he_\tx{4})({f'}_{{}_{(4')}}+{f'}_{{}_{(4')}}^{\;\;\tmu\tx{4}}\he_\tmu\he_\tx{4})(f_{{}_{(5)}}+f_{{}_{(5)}}^{\;\;\tmu\tx{5}}\he_\tmu\he_\tx{5})\dots$ and the second type $(f_{{}_{(0)}}+f_{{}_{(0)}}^{\;\;\tx{0}\ti}\he_\tx{0}\he_\ti)({f'}_{{}_{(0')}}+{f'}_{{}_{(0')}}^{\;\;\tx{0}\ti}\he_\tx{0}\he_\ti)(f_{{}_{(1)}}+f_{{}_{(1)}}^{\;\;\tx{1}\ti}\he_\tx{1}\he_\ti)\dots$ and the trivial rotation (, e.g. the rotation in the tangent space $e^{\theta_{(\tx{01})}\he_\tx{0}\he_\tx{1}}$ of the submanifold, or the rotation in the normal space $e^{\theta_{(\tx{45})}\he_\tx{4}\he_\tx{5}}$ of the submanifold), i.e. the solution of the Killing spin field equation can be written as
\begin{align*}
    \psi \;=\; & (e^{\theta_{(\tx{0}\tx{1})}\he_\tx{0}\he_\tx{1}}e^{\theta_{(\tx{0}\tx{2})}\he_\tx{0}\he_\tx{2}}\dots)(f_{{}_{(0)}}+f_{{}_{(0)}}^{\;\;\tx{0}\ti}\he_\tx{0}\he_\ti)({f'}_{{}_{(0')}}+{f'}_{{}_{(0')}}^{\;\;\tx{0}\ti}\he_\tx{0}\he_\ti)(f_{{}_{(1)}}+f_{{}_{(1)}}^{\;\;\tx{1}\ti}\he_\tx{1}\he_\ti)\dots \\
    & \qquad\qquad\qquad\quad \times (f_{{}_{(4)}}+f_{{}_{(4)}}^{\;\;\tmu\tx{4}}\he_\tmu\he_\tx{4})({f'}_{{}_{(4')}}+{f'}_{{}_{(4')}}^{\;\;\tmu\tx{4}}\he_\tmu\he_\tx{4})(f_{{}_{(5)}}+f_{{}_{(5)}}^{\;\;\tmu\tx{5}}\he_\tmu\he_\tx{5})\dots
\end{align*}

This conjecture follows from the geometric viewpoint: 
(1) There exists a neighborhood of any point in the manifold locally diffeomorphic to the Euclidean space.
(2) Any rotation operator can be expressed as the product of the basic rotations.
We can consider the local diffeomorphism of the Euclidean space as a series of deformations by the frame rotations.
Thus, any spin field $\psi$ which satisfies the Killing spin field equation locally can be written as the product of the solutions we find.

\section{Discussion}

This paper introduces a kind of spin field $\psi$, that satisfies a Killing spin field equation.
And we check that the local rotated basis $\tilde{\psi}\he_\tI\psi$ can be seen as the moving frame of a submanifold.
With the Janet-Cartan theorem, in 10-dimensional flat space this approach cam locally describe all 4-dimensional (pseudo-)Riemannian manifolds.
Through the spin field, we can see a linear relation between the connection and the extrinsic curvature of the submanifold.

At the beginning, our aim is to develop an approach using the Clifford algebra or spin field to calculate the non-linear problems in general relativity.
According to this paper, we indeed can use the spin fields to describe the geometry of space-time, but now we lack the associated physical rules to pick the appropriate spin fields.
We hope this aim can be achieved in future studies.

On the other hand, some people use the Clifford algebra to construct the unified field models \cite{Daviau:2017,Chisholm:2008,Besprosvany:2002}.\footnotemark[3]
The reasons are :
(1) the spinor can always be presented as the minimal left ideal of the spin group (which belongs to the even grade Clifford algebra);
(2) the Lie algebra in the Standard Particle Model can be embedded in certain Clifford algebras;
(3) the Clifford algebras are the more natural geometric objects than the unknown algebraic combinations.
And this paper provides another intuitive geometric interpretation that if a spinor of the elementary particle could be presented as the even grade Clifford algebras, then it should have a corresponding space-time geometry, i.e. the elementary particles could be interpreted as the nontrivial 4-dim manifolds.
We think the purely geometrical formulation of unified field theory is much closer to the idea of Einstein's unified field theory.
\footnotetext[3]{
In fact, our approach is different to theirs.
They use the Clifford algebra to represent the algebraic structure of all bosons and fermions.
But, we only use the even grade Clifford algebras to describe the fermions, and the (boson) gauge fields are not the directly represented in our Clifford algebra.
However, the spin fiedl satisfies the Killing spin field equation $\partial\psi=\dots+A\psi+\dots$ that implies the gauge field in fact related with the Clifford algebra.
We would like to point out that geometrizing the elementary particle with their results is not straightforward, but it should be worth a try.
}

\appendix

\section{Detail calculation for section \ref{sec:generating_new_solutions1}}

In the following calculation, the repeated symbolic indices obey the Einstein summation convention, except the specifying summation range with the summation symbol.
\begin{align*}
    & \psi_1 \left( \frac{1}{4}\he_\tmu\he_\tnu\omega^{{}^{(2)}\;\;\tmu\tnu}_{\;\;\;\;\lalpha} + \frac{1}{2}\he_\tmu\he_\ti H^{{}^{(2)}\;\;\ti\tmu}_{\;\;\;\;\lalpha} + \frac{1}{4}\he_\ti\he_\tj A^{{}^{(2)}\;\; \ti\tj}_{\;\;\;\;\lalpha} \right) \tilde{\psi}_1 \\
    \;=\; & \left( f + \he_\tgamma\he_\tx{5}f^{\tgamma \tx{5}} \right) \left( \frac{1}{4}\he_\tmu\he_\tnu\omega^{{}^{(2)}\;\;\tmu\tnu}_{\;\;\;\;\lalpha} + \frac{1}{2}\he_\tmu\he_\ti H^{{}^{(2)}\;\;\tmu\ti}_{\;\;\;\;\lalpha} + \frac{1}{4}\he_\ti\he_\tj A^{{}^{(2)}\;\;\ti\tj}_{\;\;\;\;\lalpha} \right) \left( f + \he_\tx{5}\he_\tdelta f^{\tdelta \tx{5}} \right) \\
    \;=\; & \frac{1}{4}\he_\tmu\he_\tnu\omega^{{}^{(2)}\;\;\tmu\tnu}_{\;\;\;\;\lalpha} \left( f \right)^2 + \frac{1}{2}\he_\tmu\he_\ti H^{{}^{(2)}\;\;\tmu\ti}_{\;\;\;\;\lalpha} \left( f \right)^2 + \frac{1}{4}\he_\ti\he_\tj A^{{}^{(2)}\;\; \ti\tj}_{\;\;\;\;\lalpha} \left( f \right)^2 \\
    & + \frac{1}{4}\he_\tmu\he_\tnu\he_\tx{5}\he_\tdelta \omega^{{}^{(2)}\;\;\tmu\tnu}_{\;\;\;\;\lalpha} \left( f\,f^{\tdelta \tx{5}} \right) + \frac{1}{2}\he_\tmu\he_\ti\he_\tx{5}\he_\tdelta H^{{}^{(2)}\;\;\tmu\ti}_{\;\;\;\;\lalpha} \left( f\,f^{\tdelta \tx{5}} \right) + \frac{1}{4}\he_\ti\he_\tj\he_\tx{5}\he_\tdelta A^{{}^{(2)}\;\; \ti\tj}_{\;\;\;\;\lalpha} \left( f\,f^{\delta 5} \right) \\
    & + \frac{1}{4}\he_\tgamma\he_\tx{5}\he_\tmu\he_\tnu\omega^{{}^{(2)}\;\;\tmu\tnu}_{\;\;\;\;\lalpha}\left( f^{\tgamma \tx{5}} f \right) + \frac{1}{2}\he_\tgamma\he_\tx{5}\he_\tmu\he_\ti H^{{}^{(2)}\;\;\tmu\ti}_{\;\;\;\;\lalpha}\left( f^{\tgamma \tx{5}} f \right) + \frac{1}{4} \he_\tgamma\he_\tx{5}\he_\ti\he_\tj A^{{}^{(2)}\;\; \ti\tj}_{\;\;\;\;\lalpha}\left( f^{\tgamma \tx{5}} f \right) \\
    & + \frac{1}{4}\he_\tgamma\he_\tmu\he_\tnu\he_\tdelta \omega^{{}^{(2)}\;\;\tmu\tnu}_{\;\;\;\;\lalpha}\left( f^{\tgamma}_{\;\;\tx{5}}f^{\tdelta \tx{5}} \right) + \frac{1}{2}\he_\tgamma\he_\tx{5}\he_\tmu\he_\ti\he_\tx{5}\he_\tdelta H^{{}^{(2)}\;\;\tmu\ti}_{\;\;\;\;\lalpha}\left( f^{\tgamma \tx{5}}f^{\tdelta \tx{5}} \right) \\
    & + \frac{1}{4}\he_\tx{5}\he_\ti\he_\tj\he_\tx{5} A^{{}^{(2)}\;\; \ti\tj}_{\;\;\;\;\lalpha}\left( f^{\tgamma \tx{5}} f_\tgamma^{\;\;\tx{5}} \right) \\
    & \textrm{\small in the following calculation, we seperate $\ti,\tj$ and $\tx{5}$, and we assume $\ti',\tj' \neq\tx{5}$} \\
    \;=\; & \frac{1}{4}\he_\tmu\he_\tnu\omega^{{}^{(2)}\;\;\tmu\tnu}_{\;\;\lalpha} \left( f \right)^2 + \frac{1}{2}\he_\tmu\he_{\ti'} H^{{}^{(2)}\;\;\tmu\ti'}_{\;\;\;\;\lalpha} \left( f \right)^2 + \frac{1}{2}\he_\tmu\he_\tx{5} H^{{}^{(2)}\;\;\tmu\tx{5}}_{\;\;\;\;\lalpha} \left( f \right)^2 + \frac{1}{4}\he_{\ti'}\he_{\tj'} A^{{}^{(2)}\;\;\ti'\tj'}_{\;\;\;\;\lalpha} \left( f \right)^2 + \frac{1}{2}\he_\tx{5}\he_{\ti'} A^{{}^{(2)}\;\;\tx{5}\ti'}_{\;\;\;\;\lalpha} \left( f \right)^2 \\
    & - \he_\tmu\he_\tx{5} \omega^{{}^{(2)}\;\;\tmu}_{\;\;\;\;\lalpha\;\;\tnu} \left( f\,f^{\tnu \tx{5}} \right) + \he_\tmu\he_\tnu H^{{}^{(2)}\;\; \tmu}_{\;\;\;\;\lalpha\;\;\tx{5}} \left( f\,f^{\tnu \tx{5}} \right) + \he_{\ti'}\he_\tx{5} H^{{}^{(2)}\;\;\;\;\ti'}_{\;\;\;\;\lalpha\tmu} \left( f\,f^{\tmu \tx{5}} \right) - \he_{\ti'}\he_\tmu A^{{}^{(2)}\;\;\;\;\ti'}_{\;\;\;\;\lalpha\tx{5}} \left( f\,f^{\tmu \tx{5}} \right) \\
    & + \frac{1}{4}\he_\tmu\he_\tnu \omega^{{}^{(2)}\;\;\tmu\tnu}_{\;\;\;\;\lalpha}\left( \sum_{\tgamma\neq \tmu,\tnu} f_{\tgamma\tx{5}} f^{\tgamma\tx{5}} \right) + \frac{1}{2}\he_\tmu\he_\tnu \omega^{{}^{(2)}\;\;\tnu}_{\;\;\;\;\lalpha\;\;\tgamma}\left( f^{\tmu}_{\;\;\tx{5}}f^{\tgamma \tx{5}} \right) - \frac{1}{2}\he_\tmu\he_\tnu \omega^{{}^{(2)}\;\;\tmu}_{\;\;\;\;\lalpha\;\;\tgamma}\left( f^{\tnu}_{\;\;\tx{5}}f^{\tgamma \tx{5}} \right) \\
    & + \frac{1}{2}\he_\tmu\he_{\ti'} H^{{}^{(2)}\;\;\tmu\ti'}_{\;\;\;\;\lalpha}\left( \sum_{\tgamma\neq\tmu} f_{\tgamma\tx{5}} f^{\tgamma\tx{5}} \right) + \he_\tmu\he_{\ti'} H^{{}^{(2)}\;\;\;\;\ti'}_{\;\;\;\;\lalpha\tgamma}\left( f^{\tgamma}_{\;\;\tx{5}}f^{\tmu \tx{5}} \right) \\
    & - \he_\tmu\he_\tx{5} H^{{}^{(2)}\;\;\tmu\tx{5}}_{\;\;\;\;\lalpha}\left( \sum_{\tgamma\neq\tmu} f_{\tgamma\tx{5}} f^{\tgamma\tx{5}} \right) - 2 \he_\tmu\he_\tx{5} H^{{}^{(2)}\;\;\;\;}_{\;\;\;\;\lalpha\tgamma\tx{5}}\left( f^{\tgamma \tx{5}}f^{\tmu \tx{5}} \right) \\
    & + \frac{1}{4}\he_{\ti'}\he_{\tj'} A^{{}^{(2)}\;\;\ti'\tj'}_{\;\;\;\;\lalpha}\left( \sum_\tgamma f_{\tgamma\tx{5}} f^{\tgamma\tx{5}} \right) - \he_\tx{5}\he_{\ti'} A^{{}^{(2)}\;\;\tx{5}\ti'}_{\;\;\;\;\lalpha}\left( \sum_\tgamma f_{\tgamma\tx{5}} f^{\tgamma\tx{5}} \right) \\
    \;=\; & \frac{1}{4}\he_\tmu\he_\tnu\left( \omega^{{}^{(2)}\;\;\tmu\tnu}_{\;\;\;\;\lalpha}\left( (f)^2 + \sum_{\tgamma\neq \tmu,\tnu} f_{\tgamma\tx{5}} f^{\tgamma\tx{5}} \right) + 2\omega^{{}^{(2)}\;\;\tnu}_{\;\;\;\;\lalpha\;\;\tgamma}(f^{\tmu \tx{5}}f^{\tgamma \tx{5}}) - 2\omega^{{}^{(2)}\;\;\tmu}_{\;\;\;\;\lalpha\;\;\tgamma}(f^{\tnu}_{\;\;\tx{5}}f^{\tgamma \tx{5}}) \right. \\
    & \qquad\quad + 2 H^{{}^{(2)}\;\;\tmu}_{\;\;\;\;\lalpha\;\;\tx{5}} \left( f\,f^{\tnu \tx{5}} \right) - 2 H^{{}^{(2)}\;\;\tnu}_{\;\;\;\;\lalpha\;\;\tx{5}} \left( f\,f^{\tmu \tx{5}} \right) \Bigg) \\
    & + \frac{1}{2}\he_\tmu\he_\tx{5}\left( H^{{}^{(2)}\;\;\tmu\tx{5}}_{\;\;\;\;\lalpha}\left( (f)^2 - 2\sum_{\tgamma\neq\tmu} f_{\tgamma\tx{5}} f^{\tgamma\tx{5}} \right) - 2 \omega^{{}^{(2)}\;\;\tmu}_{\;\;\;\;\lalpha\;\;\tgamma}(f\,f^{\tgamma \tx{5}}) + 4 H^{{}^{(2)}\;\;\;\;}_{\;\;\;\;\lalpha\tgamma\tx{5}}(f^{\tgamma \tx{5}}f^{\tmu \tx{5}}) \right) \\
    & + \frac{1}{2} \he_\tmu\he_{\ti'} \left( H^{{}^{(2)}\;\;\tmu\ti'}_{\;\;\;\;\lalpha}\left( (f)^2 + \sum_{\tgamma\neq\tmu} f_{\tgamma\tx{5}} f^{\tgamma\tx{5}} \right) + 2A^{{}^{(2)}\;\;\;\;\ti'}_{\;\;\;\;\lalpha\tx{5}} \left( f\,f^{\tmu \tx{5}} \right) + 2 H^{{}^{(2)}\;\;\;\;\ti'}_{\;\;\;\;\lalpha\tgamma}(f^{\tgamma}_{\;\;\tx{5}}f^{\tmu \tx{5}}) \right) \\
    & + \frac{1}{4}\he_{\ti'}\he_{\tj'} A^{{}^{(2)}\;\;\ti'\tj'}_{\;\;\;\;\lalpha} \\
    & + \frac{1}{2}\he_\tx{5}\he_{\ti'} \left( A^{{}^{(2)}\;\;\tx{5}\ti'}_{\;\;\;\;\lalpha} \left( (f)^2 - 2 \left( \sum_\tgamma f_{\tgamma\tx{5}} f^{\tgamma\tx{5}} \right) \right) - 2 H^{{}^{(2)}\;\;\;\;\ti'}_{\;\;\;\;\lalpha\tgamma}(f\,f^{\tgamma \tx{5}}) \right)
\end{align*}

\section{Detail calculation for section \ref{sec:generating_new_solutions2}}

In the following calculation, the repeated symbolic indices obey the Einstein summation convention, except the specifying summation range with the summation symbol.
\begin{align*}
    & \psi_1 \left( \frac{1}{4}\he_\tmu\he_\tnu\omega^{{}^{(2)}\;\;\tmu\tnu}_{\;\;\;\;\lalpha} + \frac{1}{2}\he_\tmu\he_\ti H^{{}^{(2)}\;\;\ti\tmu}_{\;\;\;\;\lalpha} + \frac{1}{4}\he_\ti\he_\tj A^{{}^{(2)}\;\; \ti\tj}_{\;\;\;\;\lalpha} \right) \tilde{\psi}_1 \\
    \;=\; & \left( f + \he_\tx{1}\he_\tk f^{\tx{1} \tk} \right) \left( \frac{1}{4}\he_\tmu\he_\tnu\omega^{{}^{(2)}\;\;\tmu\tnu}_{\;\;\;\;\lalpha} + \frac{1}{2}\he_\tmu\he_\ti H^{{}^{(2)}\;\;\tmu\ti}_{\;\;\;\;\lalpha} + \frac{1}{4}\he_\ti\he_\tj A^{{}^{(2)}\;\;\ti\tj}_{\;\;\;\;\lalpha} \right) \left( f + \he_\tm\he_\tx{1} f^{\tx{1} \tm} \right) \\
    \;=\; & \frac{1}{4}\he_\tmu\he_\tnu\omega^{{}^{(2)}\;\;\tmu\tnu}_{\;\;\;\;\lalpha} \left( f \right)^2 + \frac{1}{2}\he_\tmu\he_\ti H^{{}^{(2)}\;\;\tmu\ti}_{\;\;\;\;\lalpha} \left( f \right)^2 + \frac{1}{4}\he_\ti\he_\tj A^{{}^{(2)}\;\; \ti\tj}_{\;\;\;\;\lalpha} \left( f \right)^2 \\
    & + \frac{1}{4}\he_\tmu\he_\tnu\he_\tm\he_\tx{1} \omega^{{}^{(2)}\;\;\tmu\tnu}_{\;\;\;\;\lalpha} \left( f\,f^{\tx{1} \tm} \right) + \frac{1}{2}\he_\tmu\he_\ti\he_\tm\he_\tx{1} H^{{}^{(2)}\;\;\tmu\ti}_{\;\;\;\;\lalpha} \left( f\,f^{\tx{1} \tm} \right) + \frac{1}{4}\he_\ti\he_\tj\he_\tm\he_\tx{1} A^{{}^{(2)}\;\; \ti\tj}_{\;\;\;\;\lalpha} \left( f\,f^{\tx{1}\tm} \right) \\
    & + \frac{1}{4}\he_\tx{1}\he_\tk\he_\tmu\he_\tnu\omega^{{}^{(2)}\;\;\tmu\tnu}_{\;\;\;\;\lalpha}\left( f^{\tx{1} \tk} f \right) + \frac{1}{2}\he_\tx{1}\he_\tk \he_\tmu\he_\ti H^{{}^{(2)}\;\;\tmu\ti}_{\;\;\;\;\lalpha}\left( f^{\tx{1} \tk} f \right) + \frac{1}{4} \he_\tx{1}\he_\tk\he_\ti\he_\tj A^{{}^{(2)}\;\; \ti\tj}_{\;\;\;\;\lalpha}\left( f^{\tx{1} \tk} f \right) \\
    & + \frac{1}{4}\he_\tx{1}\he_\tmu\he_\tnu\he_\tx{1} \omega^{{}^{(2)}\;\;\tmu\tnu}_{\;\;\;\;\lalpha}\left( f^{\tx{1}}_{\;\;\tk}f^{\tx{1} \tk} \right) + \frac{1}{2}\he_\tx{1}\he_\tk\he_\tmu\he_\ti\he_\tm\he_\tx{1} H^{{}^{(2)}\;\;\tmu\ti}_{\;\;\;\;\lalpha}\left( f^{\tx{1} \tk}f^{\tx{1} \tm} \right) \\
    & + \frac{1}{4}\he_\tk\he_\ti\he_\tj\he_\tm A^{{}^{(2)}\;\; \ti\tj}_{\;\;\;\;\lalpha}\left( f_\tx{1}^{\;\;\tk} f^{\tx{1} \tm} \right) \\
    & \textrm{\small in the following calculation, we seperate $\tmu,\tnu$ and $\tx{1}$, i.e. we assume $\tmu',\tnu'\neq\tx{1}$} \\
    \;=\; & \frac{1}{4}\he_{\tmu'}\he_{\tnu'}\omega^{{}^{(2)}\;\;\tmu'\tnu'}_{\;\;\lalpha} \left( f \right)^2 + \frac{1}{2}\he_\tx{1}\he_{\tnu'}\omega^{{}^{(2)}\;\;\tx{1}\tnu'}_{\;\;\lalpha} \left( f \right)^2 + \frac{1}{2}\he_{\tmu'}\he_\ti H^{{}^{(2)}\;\;\tmu'\ti}_{\;\;\;\;\lalpha} \left( f \right)^2 + \frac{1}{2}\he_\tx{1}\he_\ti H^{{}^{(2)}\;\;\tx{1}\ti}_{\;\;\;\;\lalpha} \left( f \right)^2 + \frac{1}{4}\he_\ti\he_\tj A^{{}^{(2)}\;\;\ti\tj}_{\;\;\;\;\lalpha} \left( f \right)^2 \\
    & - \he_{\tmu'}\he_\ti \omega^{{}^{(2)}\;\;\tmu'}_{\;\;\;\;\lalpha\;\;\tx{1}} \left( f\,f^{\tx{1} \ti} \right) + \he_\ti\he_\tj H^{{}^{(2)}\;\;\;\;\ti}_{\;\;\;\;\lalpha\tx{1}} \left( f\,f^{\tx{1} \tj} \right) + \he_{\tmu'}\he_\tx{1} H^{{}^{(2)}\;\;\tmu'}_{\;\;\;\;\lalpha\;\;\tk} \left( f\,f^{\tx{1} \tk} \right) - \he_\tx{1}\he_\ti A^{{}^{(2)}\;\; \ti}_{\;\;\;\;\lalpha\;\;\tk} \left( f\,f^{\tx{1} \tk} \right) \\
    & + \frac{1}{4}\he_{\tmu'}\he_{\tnu'} \omega^{{}^{(2)}\;\;\tmu'\tnu'}_{\;\;\;\;\lalpha}\left( f_{\tx{1} \tk}f^{\tx{1} \tk} \right) + \he_\tx{1}\he_{\tmu'} \omega^{{}^{(2)}\;\;\tmu'}_{\;\;\;\;\lalpha\;\;\tx{1}}\left( f^\tx{1}_{\;\; \tk}f^{\tx{1} \tk} \right) \\
    & + \frac{1}{2}\he_{\tmu'}\he_\ti H^{{}^{(2)}\;\;\tmu'\ti}_{\;\;\;\;\lalpha}\left( \sum_{\tk\neq\ti} f_{\tx{1} \tk}f^{\tx{1} \tk} \right) + \he_\ti\he_{\tmu'} H^{{}^{(2)}\;\;\tmu'}_{\;\;\;\;\lalpha\;\;\tk}\left( f_\tx{1}^{\;\; \ti}f^{\tx{1} \tk} \right) \\
    & - \he_\tx{1}\he_\ti H^{{}^{(2)}\;\;\tx{1}\ti}_{\;\;\;\;\lalpha}\left( \sum_{\tk\neq\ti} f_{\tx{1} \tk}f^{\tx{1} \tk} \right) + 2 \he_\tx{1}\he_\ti H^{{}^{(2)}\;\;}_{\;\;\;\;\lalpha\tx{1}\tk}\left( f^{\tx{1} \tk}f^{\tx{1} \ti} \right) \\
    & + \frac{1}{4}\he_\ti\he_\tj A^{{}^{(2)}\;\;\ti\tj}_{\;\;\;\;\lalpha}\left( \sum_{\tk\neq \ti,\tj} f_{\tx{1} \tk}f^{\tx{1} \tk} \right) + \frac{1}{2}\he_\ti\he_\tj A^{{}^{(2)}\;\;\tj}_{\;\;\;\;\lalpha\;\;\tm}\left( f_\tx{1}^{\;\;\ti}f^{\tx{1}\tm} \right) - \frac{1}{2}\he_\ti\he_\tj A^{{}^{(2)}\;\;\ti}_{\;\;\;\;\lalpha\;\;\tm}\left( f_\tx{1}^{\;\;\tj}f^{\tx{1}\tm} \right) \\
    \;=\; & \frac{1}{4}\he_{\tmu'}\he_{\tnu'} \omega^{{}^{(2)}\;\;\tmu'\tnu'}_{\;\;\;\;\lalpha} \\
    & + \frac{1}{2}\he_\tx{1}\he_{\tmu'} \left( \omega^{{}^{(2)}\;\;\tx{1}\tmu'}_{\;\;\;\;\lalpha}(f)^2 - 2 \omega^{{}^{(2)}\;\;\;\;\tmu'}_{\;\;\;\;\lalpha\tx{1}}\left( f^{\tx{1}}_{\;\;\tk}f^{\tx{1}\tk} \right) - 2 H^{{}^{(2)}\;\;\tmu'}_{\;\;\;\;\lalpha\;\;\tk}(f\,f^{\tx{1}\tk}) \right) \\
    & + \frac{1}{2}\he_{\tmu'}\he_\ti\left( H^{{}^{(2)}\;\;\tmu'\ti}_{\;\;\;\;\lalpha} \left( ( f )^2 + \sum_{\tk\neq\ti} f_{\tx{1} \tk}f^{\tx{1} \tk} \right) - 2 \omega^{{}^{(2)}\;\;\tmu'}_{\;\;\;\;\lalpha\;\;\tx{1}} \left( f\,f^{\tx{1} \ti} \right) - 2 H^{{}^{(2)}\;\;\tmu'}_{\;\;\;\;\lalpha\;\;\tk}\left( f_\tx{1}^{\;\; \ti}f^{\tx{1} \tk} \right) \right) \\
    & + \frac{1}{2} \he_\tx{1}\he_\ti \left( H^{{}^{(2)}\;\;\tx{1}\ti}_{\;\;\;\;\lalpha} \left( ( f )^2 - 2 \sum_{\tk\neq\ti} f_{\tx{1} \tk}f^{\tx{1} \tk} \right) - 2 A^{{}^{(2)}\;\; \ti}_{\;\;\;\;\lalpha\;\;\tk} \left( f\,f^{\tx{1} \tk} \right) + 4 H^{{}^{(2)}\;\;}_{\;\;\;\;\lalpha\tx{1}\tk}\left( f^{\tx{1} \tk}f^{\tx{1} \ti} \right) \right) \\
    & + \frac{1}{4}\he_\ti\he_\tj \left( A^{{}^{(2)}\;\;\ti\tj}_{\;\;\;\;\lalpha} \left( ( f )^2 + \sum_{k\neq i,j} f_{\tx{1} \tk}f^{\tx{1} \tk} \right) + 2 A^{{}^{(2)}\;\;\tj}_{\;\;\;\;\lalpha\;\;\tm}\left( f_\tx{1}^{\;\;\ti}f^{\tx{1}\tm} \right) - 2 A^{{}^{(2)}\;\;\ti}_{\;\;\;\;\lalpha\;\;\tm}\left( f_\tx{1}^{\;\;\tj}f^{\tx{1}\tm} \right) \right. \\
    & \qquad + 2 H^{{}^{(2)}\;\;\;\;\ti}_{\;\;\;\;\lalpha\tx{1}} \left( f\,f^{\tx{1} \tj} \right) - 2 H^{{}^{(2)}\;\;\;\;\tj}_{\;\;\;\;\lalpha\tx{1}} \left( f\,f^{\tx{1} \ti} \right) \Bigg)
\end{align*}
%


\end{document}